\documentclass[graybox]{svmult}

\usepackage{mathptmx}       
\usepackage{helvet}         
\usepackage{courier}        
\usepackage{type1cm}        
\usepackage{makeidx}         
\usepackage{graphicx}        
\usepackage{multicol}        
\usepackage[bottom]{footmisc}

\usepackage{amsfonts}
\usepackage{enumerate}
\usepackage{amsmath}
\usepackage{amssymb}
\usepackage{amsxtra}

\makeindex             

\begin{document}

\title*{Recent progress in isoparametric functions and isoparametric hypersurfaces}
\author{Chao Qian and Zizhou Tang}
\institute{Chao Qian \at School of Mathematical Sciences, University of Chinese Academy of
Sciences, Beijing 100049, \email{qianchao@ucas.ac.cn}
\and Zizhou Tang \at The corresponding author. School of Mathematical Sciences, Laboratory of
Mathematics and Complex Systems, Beijing Normal University, Beijing 100875, \email{zztang@bnu.edu.cn}}
%
%
\maketitle

\abstract*{This paper gives a survey of recent progress in
isoparametric functions and isoparametric hypersurfaces, mainly in
two directions.
\begin{enumerate}
\item[(1)]{Isoparametric functions on Riemannian manifolds, including exotic spheres. The existences and non-existences will be considered.}
\item[(2)]{The Yau conjecture on the first eigenvalues of the embedded minimal hypersurfaces in the unit spheres. The
history and progress of the Yau conjecture on minimal isoparametric
hypersurfaces will be stated.}
\end{enumerate}
} \abstract{This paper gives a survey of recent progress in
isoparametric functions and isoparametric hypersurfaces, mainly in
two directions.
\begin{enumerate}
\item[(1)]{Isoparametric functions on Riemannian manifolds, including exotic spheres. The existences and non-existences will be considered.}
\item[(2)]{The Yau conjecture on the first eigenvalues of the embedded minimal hypersurfaces in the unit spheres. The
history and progress of the Yau conjecture on minimal isoparametric
hypersurfaces will be stated.}
\end{enumerate}
}

\section{Introduction}
\label{sec:1} E. Cartan was the pioneer who made a comprehensive
study of isoparametric functions ( hypersurfaces ) on the unit
spheres. In the past decades, the study of isoparametric functions (
hypersurfaces ) has became a highly influential field in
differential geometry. For a systematic and complete survey of
isoparametric functions ( hypersurfaces ) and their generalizations,
we recommend \cite{Th00}, \cite{Ce08} and \cite{Ch12}. Very recently,
Cecil-Chi-Jensen, Immervoll and Chi obtained classification results
for isoparametric hypersurfaces with four distinct principal
curvatures in the unit spheres, except for one case ( c.f.
\cite{CCJ07}, \cite{Imm08} and \cite{Ch13} ). As for that with six
distinct principal curvatures, Miyaoka showed the homogeneity and
hence the classification ( c.f. \cite{Miy13} ).

The note is organized as follows. In Section \ref{sec:2}, we first
recall some basic notations and fundamental theory of isoparametric
functions on Riemannian manifolds. Next we introduce exotic spheres
and investigate the existences and non-existences of isoparametric
functions on exotic spheres. Section \ref{sec:3} will be concerned
with the progress of the well known Yau conjecture on the first
eigenvalues of embedded minimal hypersurfaces in the unit spheres,
especially on the minimal isoparametric case ( being isoparametric
implies embedding ). Moreover, the first eigenvalues of the focal
submanifolds are also taken into account. In the end, related topics
and applications are described in Section \ref{sec:4}.

\section{Exotic spheres and isoparametric functions}
\label{sec:2}

We start with definitions. Let $N$ be a connected complete
Riemannian manifold. A non-constant smooth function $f: N\rightarrow
\mathbb{R}$ is called \emph{transnormal} if there is a smooth
function $b: \mathbb{R}\rightarrow\mathbb{R}$ such that
\begin{equation}\label{iso1}
|\nabla f|^2=b(f),
\end{equation}
where $\nabla f$ is the gradient of $f$. If moreover there is a continuous function
$a: \mathbb{R}\rightarrow\mathbb{R}$ such that
\begin{equation}\label{iso2}
\triangle f=a(f),
\end{equation}
where $\triangle f$ is the Laplacian of $f$, then $f$ is called
\emph{isoparametric} ( cf. \cite{Wa87} ). Each regular level
hypersurface is called an \emph{isoparametric hypersurface}. The two
equations of the function $f$ mean that regular level hypersurfaces
are parallel and have constant mean curvatures. According to Wang
\cite{Wa87}, a transnormal function $f$ on a complete Riemannian
manifold has no critical value in the interior of $\mathrm{Im}f$.
The preimage of the maximum ( resp. minimum ), if it exists, of an
isoparametric ( or transnormal ) function $f$ is called the
\emph{focal set} of $f$, denoted by $M_{+}$ ( resp. $M_{-}$ ).

Since the work of Cartan ( \cite{Ca38}, \cite{Ca39} ) and
M{\"u}nzner ( \cite{Mu80} ), the subject of isoparametric
hypersurfaces in the unit spheres is rather fascinating to
geometers. We refer to \cite{CR85} for the development of this
subject. Up to now, the classification has almost been completed as
mentioned in Section \ref{sec:1}.

In general Riemannian manifolds, the classification problem is far
from being touched. Wang \cite{Wa87} firstly took up a systematic
study of isoparametric functions on general Riemannian manifolds,
and similar to the case in a unit sphere, proved or claimed a series
of beautiful results. The structural result for transnormal
functions is stated as follows.
\begin{theorem}\label{tube}(\cite{Wa87})
Let $N$ be a connected complete Riemannian manifold and $f$ a transnormal function on $N$. Then
\begin{itemize}
\item{The focal sets of $f$ are smooth submanifolds ( may be disconnected ) of $N$;}
\item{Each regular level set of $f$ is a tube over either of the focal sets ( the
          dimensions of the fibers may differ on different connected components ).}
\end{itemize}
\end{theorem}

The above theorem shows that the existence of a transnormal function
on a Riemannian manifold $N$ restricts strongly its topology.

In the first part of \cite{GT13}, Ge and Tang improved the
fundamental theory of isoparametric functions on Riemannian
manifolds. Given a transnormal function $f: N\rightarrow
\mathbb{R}$, we denote by $C_1(f)$ the set where $f$ attains its
global maximum value or global minimum value, by $C_2(f)$ the union
of singular level sets of $f$, \emph{i.e.}, $C_2(f)=\{p\in N| \nabla
f(p)=0\}$, and for any regular value $t$ of $f$, by $C^t_3(f)$ the
focal set of the level hypersurface $M_t:=f^{-1}(t)$, \emph{i.e.},
the set of singular values of the normal exponential map. From
\cite{Wa87}, it follows that $C_1(f)=C_2(f)=M_{-}\cup M_{+}$, and
$C^{t_1}_3(f)=C^{t_2}_3(f)$ for any two regular level hypersurfaces
which will be thus denoted simply by $C_3(f)$. Moreover, one can see
that $C_3(f)\subset C_1(f)=C_2(f)$. Then Ge and Tang proved
\begin{theorem}\label{proper equiv-def}(\cite{GT13})
Each component of $M_{\pm}$ has codimension not less than $2$ if and
only if~$C_3(f)=C_1(f)=C_2(f)$. Moreover in this case, each level
set $M_t$ is connected. If in addition $N$ is closed and $f$ is
isoparametric, then at least one isoparametric hypersurface is
minimal in $N$.
\end{theorem}

Indeed, there exists example of an isoparametric funtion $f$
satisfying $C_3(f)\varsubsetneq C_1(f)=C_2(f)$ ( c.f. \cite{GT13} ).
For this case, the focal sets of the isoparametric function are not
really focal sets of the level hypersurface. Hence, in \cite{GT13},
a transnormal (isoparametric) function $f$ is called \emph{proper}
if the focal sets have codimension not less than $2$. It seems that
a properly transnormal (isoparametric) function is exactly what we
should concern in geometry. Furthermore, in \cite{GT13}, they
observed three elegant ways to construct examples of isoparametric
functions, \emph{i.e.},
\begin{itemize}
\item{ For a Riemannian manifold $(N, ds^2)$ with an isoparametric
function $f$, take a special conformal deformation
$\widetilde{ds^2}=e^{2u(f)}ds^2$. Then $f$ is also isoparametric on
$(N, \widetilde{ds^2})$;}
\item{ For a cohomogeneity one manifold
$(N,G)$ with a $G$-invariant metric, one can get isoparametric
functions on $N$;}
\item{ For a Riemannian submersion $\pi:
E\rightarrow B$ with minimal fibers, if $f$ is an isoparametric
function on $B$, then so is $F:=f\circ \pi$ on $E$.}
\end{itemize}

Applying these methods, interesting results and abundant examples
are acquired, especially, isoparametric functions on Brieskorn
varieties and on isoparametric hypersurfaces of spheres are
obtained.

As a continuation of \cite{GT13}, they made new contributions in
\cite{GT14}. First, for a properly isoparametric function, they
proved that at least one isoparametric hypersurface is minimal if
the ambient space $N$ is closed in Theorem \ref{proper equiv-def}.
By using the Riccati equation, they can further show that such a
minimal isoparametric hypersurface is also unique if $N$ has
positive Ricci curvature. Next, by expressing the shape operator
$S(t)$ of $M_t$ as a power series, they gave a complete proof to
Theorem D of \cite{Wa87} ( no proof there; compare with \cite{Ni97}
and \cite{Miy13a} ).
\begin{theorem}\label{minimal focal}(\cite{GT14})
The focal sets $M_{\pm}$ of an isoparametric function $f$ on a
complete Riemannian manifold $N$ are minimal submanifolds.
\end{theorem}

Meanwhile, Ge and Tang also established the following theorem, which
is a generalization of the spherical case to general Riemannian
manifolds.
\begin{theorem}\label{austere focal}(\cite{GT14})
Suppose that each isoparametric hypersurface $M_t$ has constant
principal curvatures with respect to the unit normal vector field in
the direction of $\nabla f$. Then each of the focal sets $M_{\pm}$
has common constant principal curvatures in all normal directions,
\emph{i.e.}, the eigenvalues of the shape operator are constant and
independent of the choices of the point and unit normal vector of
$M_{\pm}$.
\end{theorem}

Owning to the rich and beautiful topological and geometric
properties of isoparametric functions on Riemannian manifolds, Ge
and Tang initiated the study of isoparametric functions on exotic
spheres in \cite{GT13}.

Recall that an $n$-dimensional smooth manifold $\Sigma^n$ is called
an \emph{exotic n-sphere} if it is homeomorphic but not
diffeomorphic to $S^n$. It is J. Milnor \cite{Mil56} who firstly
discovered an exotic 7-sphere which is an $S^3$-bundle over $S^4$.
Later, Kervaire and Milnor \cite{KM63} computed the group of
homotopy spheres in each dimension greater than four which implies
that there exist exotic spheres in infinitely many dimensions and in
each dimension there are at most finitely many exotic spheres. In
particular, ignoring orientation there exist 14 exotic 7-spheres, 10
of which can be exhibited as $S^3$-bundles over $S^4$, the so-called
\emph{Milnor spheres}. However, in dimension four, the question of
whether an exotic 4-sphere exists remains open, which is the so
called \emph{smooth Poincar\'{e} conjecture} ( c.f. \cite{JW08} ).

Since the discovery of exotic spheres by Milnor, a very intriguing
problem is to interpret the geometry of them ( c.f. \cite{GM74},
\cite{GZ00}, \cite{BGK05}, and \cite{JW08} ). In \cite{GT13}, using
isoparametric ( even transnormal ) functions to attack the smooth
Poincar\'{e} conjecture in dimension four, Ge and Tang showed the
following theorem.
\begin{theorem}\label{homotopy 4-sphere}(\cite{GT13})
Suppose $\Sigma^4$ is a homotopy $4$-sphere and it admits a transnormal
function under some metric. Then $\Sigma^4$ is diffeomorphic to $S^4$.
\end{theorem}

Note that a \emph{homotopy n-sphere} is a smooth manifold with the
same homotopy type as $S^n$. Freedman \cite{Fr82} showed that any
homotopy 4-sphere is homeomorphic to $S^4$. As a result of this, the
above Theorem \ref{homotopy 4-sphere} says equivalently that there
exists no transnormal function on any exotic 4-sphere if it exists.
In contrast to the non-existence result in dimension four, Ge and
Tang also constructed many examples of isoparametric functions on
the Milnor spheres. Furthermore, by projecting an $S^3$-invariant
isoparametric function on the symplectic group $Sp(2)$ with a
certain left invariant metric, they constructed explicitly a
properly transnormal but not an isoparametric function on the
Gromoll-Meyer sphere with two points as the focal sets. Inspired by
this example, they posed a question that whether there is an
isoparametric function on the Gromoll-Meyer sphere or on any exotic
n-sphere ($n>4$) with two points as the focal sets. More generally,
they posed the following:
\begin{problem}\label{problem}(\cite{GT13})
Does there always exist a properly isoparametric function on an
exotic sphere $\Sigma^n$ ($n>4$) with the focal sets being
those occurring on $S^n$?
\end{problem}

To answer the Problem \ref{problem}, Qian and Tang developed a
general way to construct metrics and isoparametric functions on a
given manifold in \cite{QT13}, which is based on a simple and useful
observation that a transnormal function on a complete Riemannian
manifold is necessarily a Morse-Bott function ( c.f. \cite{Wa87} ).
As is well known, a Morse-Bott function is a generalization of a
Morse function, and it admits critical submanifolds satisfying a
certain non-degenerate condition on normal bundles. In \cite{QT13},
the following fundamental construction is given, whose proof depends
heavily on Moser's volume element theorem.
\begin{theorem}\label{Main}(\cite{QT13})
Let $N$ be a closed connected smooth manifold and $f$ a Morse-Bott
function on $N$ with the critical set $C(f)=M_+\sqcup M_-$, where
$M_+$ and $M_-$ are both closed connected submanifolds of
codimensions more than 1. Then there exists a metric on $N$ so that
$f$ is an isoparametric function.
\end{theorem}

It follows from a theorem of S. Smale that

\begin{corollary}\label{homosphere}(\cite{QT13})
Every homotopy $n$-sphere with $n>4$ admits a metric and an
isoparametric function with 2 points as the focal sets.
\end{corollary}
\begin{remark}
Corollary \ref{homosphere} answers partially the above Problem \ref{problem}.
\end{remark}

Moreover, metrics and isoparametric functions on homotopy spheres
and on the Eells-Kuiper projective planes can also be constructed so
that at least one component of the critical set is not a single
point.

In addition to the above existence theorem on homotopy spheres, on the other side,
the following non-existence results was also proved.
\begin{theorem}\label{nonisopara}(\cite{QT13})
Every odd dimensional exotic sphere admits no totally isoparametric
functions with 2 points as the focal set.
\end{theorem}
Recall that a \emph{totally isoparametric} function is an
isoparametric function so that each regular level hypersurface has
constant principal curvatures, as defined in \cite{GTY11}. As it is
well known, an isoparametric function on a unit sphere must be
totally isoparametric.
\begin{remark}
According to \cite{HH67} and \cite{St96}, there exists at leat one exotic
Kervaire sphere $\Sigma^{4m+1}$ which has a cohomogeneity one action.
Consequently, $\Sigma^{4m+1}$ admits a
totally isoparametric function $f$ under an invariant metric ( c.f. \cite{GT13} ).
However, each component of the focal set of $f$ is not just a point,
but a smooth submanifold. Hence, the assumption on the focal set in
Theorem \ref{nonisopara} is essential.
\end{remark}

In light of the above Theorem \ref{nonisopara}, it is reasonable to
ask
\begin{problem}
Does there exist an even dimensional exotic sphere
$\Sigma^{2n}(n>2)$ admitting a metric and a totally isoparametric
function with 2 points as the focal set?
\end{problem}

In the last section of \cite{QT13}, both existence and non-existence
results of isoparametric functions on some homopoty spheres which
also have $SC^p$-property were investigated. A Riemannian manifold
has $SC^{p}$-property if every geodesic issuing from the point $p$
is closed and has the same length ( \cite{Be78} ). For some even
dimensional homotopy spheres, the existence theorem in \cite{QT13}
improves a beautiful result of B\'{e}rard-Bergery \cite{BB77}.

Recently, Tang and Zhang \cite{TZ14} solved a problem of
B\'{e}rard-Bergery and Besse. That is, they showed that every
Eells-Kuiper quaternionic projective plane carries a Riemannian
metric with $SC^p$-property for a certain point $p$. Thus, it is
interesting to know whether there is a metric on every Eells-Kuiper
quaternionic projective plane which not only has the
$SC^p$-property, but also admits a certain isoparametric function (
c.f. \cite{QT13} ).

\section{Yau conjecture on the first eigenvalue and isoparametric foliations}
\label{sec:3}
The Laplace-Beltrami operator is one of the most important operators acting on $C^{\infty}$ functions
on a Riemannian manifold. Over several decades,
research on the spectrum of the Laplace-Beltrami operator has always been a
core issue in the study of geometry.  For instance, the geometry of closed minimal submanifolds in
the unit sphere is closely related to the eigenvalue problem.

Let $(M^n,g)$ be an $n$-dimensional compact connected Riemannian manifold without boundary
and $\Delta $ be the Laplace-Beltrami operator acting on a $C^{\infty}$ function $f$ on
$M$ by $\Delta f$ $=-$ div$(\nabla f)$,
the negative of divergence of the gradient $\nabla f$.
It is well known that $\Delta$ is an elliptic operator and has a discrete spectrum
$$\{0=\lambda_0(M)<\lambda_1(M)\leq \lambda_2(M)\leq \cdots \leq\lambda_k(M),\cdots, \uparrow \infty\}$$
with each eigenvalue repeated a number of times equal to its
multiplicity. As usual, we call $\lambda_1(M)$ the first eigenvalue
of $M$. When $M^n$ is a minimal hypersurface in the unit sphere
$S^{n+1}(1)$, it follows from Takahashi Theorem that $\lambda_1(M)$
is not greater than $n$. Consequently, S.T. Yau posed in 1982 the
following conjecture: \vspace{1mm}

\noindent \textbf{Yau conjecture} (\cite{Ya82}):\,\, {\itshape The
first eigenvalue of every closed embedded minimal hypersurface $M^n$
in the unit sphere $S^{n+1}(1)$ is just $n$. } \vspace{1mm}

In 1983, Choi and Wang made the most significant breakthrough to
this conjecture ( \cite{CW83} ). To be precise, they showed that the
first eigenvalue of every (embedded) closed minimal hypersurface in
$S^{n+1}(1)$ is not smaller than $n\over 2$. Usually, the
calculation of the spectrum of the Laplace-Beltrami operator, even
of the first eigenvalue, is rather complicated and difficult. Up to
now, the Yau conjecture is far from being solved even in dimension
two.

It was proved in ( \cite{MR86} ) that if the Yau conjecture is true
for the torus of dimension two, then the Lawson conjecture holds,
that is to say, the only minimally embedded torus in $S^3(1)$ is the
Clifford torus. In fact, the Lawson conjecture has been a
challenging problem for more than 40 years, and recently it was
solved by S. Brendle ( c.f. \cite{Br13} ).

In this note, we pay attention to a little more restricted problem
of the Yau conjecture for closed minimal isoparametric hypersurfaces
$M^n$ in $S^{n+1}(1)$.

Recall that an isoparametric hypersurface $M^n$ in the unit sphere
$S^{n+1}(1)$ must have constant principal curvatures ( c.f.
\cite{Ca38}, \cite{Ca39}, \cite{CR85} ). Let $\xi$ be a unit normal
vector field along $M^n$ in $S^{n+1}(1)$, $g$ the number of distinct
principal curvatures of $M$, $\cot \theta_{\alpha}~
(\alpha=1,...,g;~ 0<\theta_1<\cdots<\theta_{g} <\pi)$ the principal
curvatures with respect to $\xi$ and $m_{\alpha}$ the multiplicity
of $\cot \theta_{\alpha}$. Using a brilliant topological method,
M\"{u}nzner ( c.f. \cite{Mu80} ) proved the remarkable result that
the number $g$ must be $1, 2, 3, 4$ or $6$;
$m_{\alpha}=m_{\alpha+2}$ (indices mod $g$);
$\theta_{\alpha}=\theta_1+\frac{\alpha-1}{g}\pi$ $(\alpha = 1,...,
g)$ and when $g$ is odd, $m_1=m_2$.

In order to attack the Yau conjecture, Muto-Ohnita-Urakawa (
\cite{MOU84} ) and Kotani ( \cite{Ko85} ) made a breakthrough for
some minimal homogeneous ( automatically isoparametric )
hypersurfaces. More precisely, they verified the Yau conjecture for
all minimal homogeneous hypersurfaces with $g= 1, 2, 3, 6$. However,
when it came to the case $g=4$, they were only able to deal with the
cases $(m_1, m_2)=(2,2)$ and $(1, k)$.

Furthermore, Muto ( \cite{Mu88} ) proved that the Yau conjecture is
also true for some families of minimal inhomogeneous isoparametric
hypersurfaces with $g=4$. This remarkable result contains many
inhomogeneous isoparametric hypersurfaces. However, there is no
result in \cite{Mu88} for isoparametric hypersurfaces with
$\min(m_1, m_2)> 10$.

Based on all results mentioned above and the classification of
isoparametric hypersurfaces in $S^{n+1}(1)$ ( c.f. \cite{CCJ07},
\cite{Imm08}, \cite{Ch13}, \cite{DN85} and \cite{Miy13} ), Tang and
Yan \cite{TY13} completely solved the Yau conjecture on the minimal
isoparametric case by establishing the following
\begin{theorem}\label{thm1 hypersurface}(\cite{TY13})
\emph{Let $M^n$ be a closed minimal isoparametric hypersurface in
the unit sphere $S^{n+1}(1)$ with $g=4$ and $m_1, m_2\geq 2$. Then
$$\lambda_1(M^n)=n.$$ }
\end{theorem}

\begin{remark}
Theorem \ref{thm1 hypersurface} depends only on the values of $(m_1,
m_2)$. In particular, it covers the unclassified case $g=4$, $(m_1,
m_2)=(7, 8)$.
\end{remark}

\begin{remark}\label{Chern conjecture}
A purported conjecture of Chern states that a closed, minimally
immersed hypersurface in $S^{n+1}(1)$, whose second fundamental form
has constant length, is isoparametric ( c.f. \cite{GT12} ). If this
conjecture is proven, Theorem \ref{thm1 hypersurface} would have
settled the Yau conjecture for the minimal hypersurface whose second
fundamental form has constant length, which gives more confidence in
the Yau conjecture.
\end{remark}

Indeed, the more fascinating part of \cite{TY13} was to determine
the first eigenvalues of the focal submanifolds in $S^{n+1}(1)$. To
state their result clearly, let us recall some notations. Given an
isoparametric hypersurface $M^n$ in $S^{n+1}(1)$ and a smooth field
$\xi$ of unit normals to $M$, for each $x\in M$ and $\theta\in
\mathbb{R}$, we can define $\phi_{\theta}: M^n\rightarrow
S^{n+1}(1)$ by
$$\phi_{\theta}(x)=\cos \theta~ x +\sin \theta~ \xi(x).$$
Clearly, $\phi_{\theta}(x)$ is the point at an oriented distance
$\theta$ to $M$ along the normal geodesic through $x$. If
$\theta\neq \theta_{\alpha}$ for any $\alpha=1,...,g$,
$\phi_{\theta}$ is a parallel hypersurface to $M$ at an oriented
distance $\theta$, which we will denote by $M_{\theta}$
henceforward. If $\theta= \theta_{\alpha}$ for some
$\alpha=1,...,g$, it is easy to find that for any vector $X$ in the
principal distributions $E_{\alpha}(x)=\{X\in T_xM ~|~A_{\xi}X=\cot
\theta_{\alpha}X\}$, where $A_{\xi}$ is the shape operator with
respect to $\xi$, $(\phi_{\theta})_{\ast}X=0$. In other words, in
case $\cot \theta=\cot \theta_{\alpha}$ is a principal curvature of
$M$, $\phi_{\theta}$ is not an immersion, which is actually a
\emph{focal submanifold} of codimension $m_{\alpha}+1$ in
$S^{n+1}(1)$.

M\"{u}nzner asserted that there are only two distinct focal
submanifolds in a parallel family of isoparametric hypersurfaces,
regardless of the number of distinct principal curvatures of $M$;
and every isoparametric hypersurface is a tube of constant radius
over each focal submanifold. Denote by $M_1$ the focal submanifold
in $S^{n+1}(1)$ at an oriented distance $\theta_1$ along $\xi$ from
$M$ with codimension $m_1+1$, $M_2$ the focal submanifold in
$S^{n+1}(1)$ at an oriented distance $\frac{\pi}{g}-\theta_1$ along
$-\xi$ from $M$ with codimension $m_2+1$. In virtue of Cartan's
identity, one sees that the focal submanifolds $M_1$ and $M_2$ are
both minimal in $S^{n+1}(1)$ ( c.f. \cite{CR85} ).

Another main result of \cite{TY13} concerning the first eigenvalues
of focal submanifolds in the non-stable range, is now stated as
follows.

\begin{theorem}\label{thm2 focal submanifold}(\cite{TY13})
\emph{Let $M_1$ be the focal submanifold of an isoparametric
hypersurface with $g=4$ in $S^{n+1}(1)$ with codimension $m_1+1$. If
$\dim M_1\geq \frac{2}{3}n+1$, then
$$\lambda_1(M_1)=\dim M_1$$ with multiplicity $n+2$.
A similar conclusion holds for $M_2$ under an analogous condition.}
\end{theorem}

We emphasize that the assumption $\dim M_1\geq \frac{2}{3}n+1$ in
Theorem \ref{thm2 focal submanifold} is essential. For instance,
Solomon \cite{So92} constructed an eigenfunction on the specific
focal submanifolds $M_2$ of OT-FKM-type ( we will explain it
immediately ), which has $4m$ as an eigenvalue. In some case, $4m$
is less than the dimension of $M_2$.

As an example, Theorem \ref{thm2 focal submanifold} implies that
each focal submanifold of isoparametric hypersurfaces with $g=4$,
$(m_1, m_2)=(7,8)$ has its dimension as the first eigenvalue.

We need to recall the construction of the isoparametric
hypersurfaces of OT-FKM-type. For a symmetric Clifford system
$\{P_0,\cdots,P_m\}$ on $\mathbb{R}^{2l}$, \emph{i.e.}, $P_i$'s are
symmetric matrices satisfying $P_iP_j+P_jP_i=2\delta_{ij}I_{2l}$,
Ferus, Karcher and M\"{u}nzner (\cite{FKM81}) constructed a
polynomial $F$ on $\mathbb{R}^{2l}$:
\begin{eqnarray*}\label{FKM isop. poly.}
&&\qquad F:\quad \mathbb{R}^{2l}\rightarrow \mathbb{R}\nonumber\\
&&F(x) = |x|^4 - 2\displaystyle\sum_{i = 0}^{m}{\langle
P_ix,x\rangle^2}.
\end{eqnarray*}
For $f=F|_{S^{2l-1}}$, define $M_1=f^{-1}(1)$, $M_2=f^{-1}(-1)$,
which have codimensions $m+1$ and $l-m$ in $S^{n+1}(1)$,
respectively.

For focal submanifold $M_1$ of OT-FKM-type, the only unsettled
multiplicities in \cite{TY13} are $(m_1, m_2)=(1, 1), (4, 3), (5,
2)$. And for the $(4,3)$ case, there exist only one homogeneous and
one inhomogeneous examples.

Finally, Tang and Yan \cite{TY13} proposed the following problem,
which could be regarded as an extension of the Yau conjecture.
\vspace{2mm}
\begin{problem}\label{Generalized Yau conjecture}(\cite{TY13}) Let $M^d$ be a closed embedded
minimal submanifold in the unit sphere $S^{n+1}(1)$. If the
dimension $d$ of $M^d$ satisfies $d\geq \frac{2}{3}n+1$, then
$$\lambda_1(M^d)=d.$$
\end{problem}

Later, Tang, Xie and Yan \cite{TXY14} took chance to solve the
unsolved cases in \cite{TY13} and considered the case with g=6.
First, by applying the similiar method as in \cite{TY13}, they got
the following theorem for the case with $g=6$ which contains more
information than that in \cite{MOU84} and does not depend on the
classification result of Miyaoka ( \cite{Miy13} ).

\begin{theorem}\label{prop 1 of TXY}(\cite{TXY14})
Let $M^{12}$ be a closed minimal isoparametric hypersurface in
$S^{13}(1)$ with $g=6$ and $(m_1, m_2)=(2, 2)$. Then
$$\lambda_1(M^{12})=12$$
with multiplicity $14$. Furthermore, the following inequality holds
$$\lambda_k(M^{12})>\frac{3}{7}~\lambda_k(S^{13}(1)), \qquad k=1, 2, \cdots.$$
\end{theorem}

And for focal submanifolds $M_1$ of OT-FKM-Type, they solved two left cases
and proved
\begin{theorem}\label{thm2 TXY focal g=4}(\cite{TXY14})
For the focal submanifold $M_1$ of OT-FKM-type in $S^5(1)$ with
$(m_1, m_2)=(1, 1)$,
$$\lambda_1(M_1)=\dim M_1=3$$
with multiplicity $6$; for the focal submanifold $M_1$ of
homogeneous OT-FKM-type in $S^{15}(1)$ with $(m_1, m_2)=(4, 3)$,
$$\lambda_1(M_1)=\dim M_1=10$$
with multiplicity $16$.
\end{theorem}

At last, in the case with $g=6$, by a deep investigation into the
shape operator of the focal submanifolds, they obtained estimates on
the first eigenvalue. Particularly, for one of the focal
submanifolds with $g=6$, $m_1=m_2=2$, the first eigenvalue is equal
to its dimension. It gives an affirmative answer to Problem
\ref{Generalized Yau conjecture} in this case.
\section{Related topics and applications}
\label{sec:4}

The connection between geometry of Riemannian manifolds with
positive scalar curvatures and surgery theory is quite a deep
subject which has attracted widely attention. The most important
aspect of this field is the original discovery of Gromov-Lawson and
of Schoen-Yau.

Motivated by the Schoen-Yau-Gromov-Lawson surgery theory on metrics
of positive scalar curvature, Tang, Xie and Yan \cite{TXY12}
constructed a double manifold associated with a minimal
isoparametric hypersurface in the unit sphere. The resulting double
manifold carries a metric of positive scalar curvature and an
isoparametric foliation as well. To investigate the topology of the
double manifolds, they used topological K-theory and the
representation of the Clifford algebra for the OT-FKM-type, and
determined completely the isotropy subgroups of singular orbits for
homogeneous case. Here we note that, as it is well known, a
homogeneous ( isoparametric ) hypersurface in the unit sphere can be
characterized as a principal orbit of isotropy representation of
some symmetric space of rank two.

In the last part of this section, we describe an application of
isoparametric foliation to Willmore submanifolds. By definition, a
Willmore submanifold ( in the unit sphere ) is the critical point of
the Willmore functional. In particular, every minimal surface in the
unit sphere is automatically Willmore; in other words, Willmore
surfaces are a generalization of minimal surfaces in the unit
sphere. However, examples of Willmore submanifolds in the unit
sphere are rare in the literature.

Qian, Tang and Yan ( \cite{TY12}, \cite{QTY13} ) proved that each
focal submanifold of isoparametric hypersurface ( not only
OT-FKM-type ) in the unit sphere with $g=4$ is a Willmore
submanifold. For $g=1, 2, 3$, the conclusion above is clearly valid.
As for $g=6$, the conclusion should be also true.

Recall that the focal submanifolds are minimal in unit spheres. It
is worth noting that an Einstein manifold minimally immersed in the
unit sphere is a Willmore submanifold. A natural problem arises:
whether the focal submanifolds are Einstein? To this problem with
$g=4$, \cite{TY12} and \cite{QTY13} gave a complete answer,
depending on the classification results. In other words, they dealt
with this problem in two cases--homogeneous type and OT-FKM-type.
%
%
%

\end{document}